\documentclass[12pt]{article}

\usepackage[applemac]{inputenc}
\usepackage{amsmath,amssymb,fullpage}
\usepackage{showlabels}
\usepackage{color}

\newcommand{\dproof}{\noindent {Proof.} \quad}
\newcommand{\fproof}{\hfill $\square$ \bigskip}

\newtheorem{definition}{Definition}[section]
\newtheorem{example}{Example}[section]
\newtheorem{theorem}[definition]{Theorem}
\newtheorem{problem}[definition]{Problem}
\newtheorem{remark}[definition]{ \it Remark}
\newtheorem{coro}[definition]{Corollary}
\newtheorem{proposition}[definition]{Proposition}

\numberwithin{equation}{section}

\def\1B{\text{1\!\!I}}

\begin{document}
\date{18 May 2016}
\title{ Optimal insider control and semimartingale decompositions under enlargement of filtration}

\author{
Olfa Draouil$^{1}$ and Bernt \O ksendal$^{2,3}$}

\footnotetext[1]{Department of Mathematics, University of Tunis El Manar, Tunis, Tunisia.
Email: {\tt olfadraouil@hotmail.fr}}

\footnotetext[2]{Department of Mathematics, University of Oslo, P.O. Box 1053 Blindern, N--0316 Oslo, Norway. Email: {\tt oksendal@math.uio.no}}

\footnotetext[3]{This research was carried out with support of CAS - Centre for Advanced Study, at the Norwegian Academy of Science and Letters, within the research program SEFE, and with support of the Norwegian Research Council, within the research project Challenges in Stochastic Control, Information and Applications (STOCONINF), project number 250768/F20.}
\maketitle
\begin{abstract}
We combine stochastic control methods, white noise analysis and Hida-Malliavin calculus applied to the Donsker delta functional to obtain explicit representations of semimartingale decompositions under enlargement of filtrations. Some of the expressions are more explicit than previously known. The results are illustrated by examples.
\end{abstract}



\section{Introduction}
The purpose of this paper is twofold:
\begin{itemize}
\item
We introduce a new approach to enlargement of filtration problems, based on combining several optimal control methods.
\item
We show that this approach can in some cases give more explicit results than known before.
\end{itemize}

The system we consider, is described by a stochastic differential equation driven by a Brownian motion $B(t)$ and an independent compensated Poisson random measure $\tilde{N}(dt,d\zeta)=N(dt,d\zeta)-\nu(d\zeta)dt$, $\nu$ being the L\' evy measure of the Poisson random measure $N$. The processes are jointly defined on a filtered probability space $(\Omega, \mathbb{F}=\{ \mathcal{F}_t \}_{t \geq 0},\mathbf{P})$ satisfying the usual conditions where $\Omega=\mathcal{S}'(\mathbb{R})$ and $\mathbf{P}$ is the Gaussian measure on $\mathcal{S}'(\mathbb{R})$. Here, and throughout the paper, $\mathcal{F}_t$ denotes the sigma-algebra generated by $\{B(s)\}_{s\leq t}$ and $\{N(s, \cdot)\}_{s\leq t}$. Throughout this paper we assume that the inside information is of \emph{initial enlargement} type. Specifically, we assume  that the inside filtration $\mathbb{H}$ has the form

\begin{equation}\label{eq1.1}
 \mathbb{H}= \mathbb{H}_Y= \{ \mathcal{H}_t\}_{t \geq 0}, \text{ where } \mathcal{H}_t = \mathcal{F}_t \vee \sigma(Y)
\end{equation}
for all $t$, where $Y$ is a given $\mathcal{F}_{T_0}$-measurable random variable, for some $T_0 > 0$ (constant).
In order to satisfy the usual hypotheses we redefine $\mathcal{H}_{t}= \mathcal{H}_{t^+}=\bigcap_{s>t}\mathcal{H}_s$.
We also assume that the Donsker delta functional of $Y$, $\delta_Y(y)$, exists as an element of the Hida space $(\mathcal{S})^{\ast}$ of stochastic distributions (see Section 2) and that
\begin{align}\label{eq1.2a}
&\mathbb{E}[\delta_Y(.)|\mathcal{F}_t]\in\mathbf{L}^2(m\times\mathbf{P}) \text{ and } \nonumber\\
&\mathbb{E}[D_t\delta_Y(.)|\mathcal{F}_t]\in\mathbf{L}^2(m\times\mathbf{P})\text{ and }\nonumber\\
& \mathbb{E}[D_{t,z}\delta_Y(.)|\mathcal{F}_t]\in\mathbf{L}^2(m\times\nu\times\mathbf{P}),
\end{align}
where $D_t$ and $D_{t,\zeta}$ denote the Hida-Malliavin derivatives with respect to $B(\cdot)$ and $\tilde{N}(\cdot,\cdot)$, respectively. We refer to \cite{DO1} for more information about $(\mathcal{S})^*$, white noise theory and Hida-Malliavin derivatives. As pointed out in \cite{DO1}, conditions \eqref{eq1.2a} hold automatically if $Y$ has the form
\begin{equation}\label{eq1.3}
 Y= \int_0^{T_0} \varphi(s)dB(s) + \int_0^{T_0} \psi(s,\zeta) \tilde{N}(ds,d\zeta),
 \end{equation}
for some deterministic functions  $\varphi, \psi$ such that  $\varphi \in L^2[m]$ and $\psi \in L^2[m \times \nu]$, where $m$ is Lebesgue measure on $[0,T_0]$,
and we assume that the L\'  evy measure $\nu$ satisfies the condition
\begin{equation}\label{eq1.3a}
\int_{\mathbb{R}}\zeta^2d\nu(\zeta)<\infty.
\end{equation}
The condition \eqref{eq1.3a} is not needed if $\int_0^{T_0}\varphi(t)^2dt > 0$.


The semimartingale decomposition problem we consider, is the following:
\begin{problem}
Find an $\mathbb{H}$-adapted process $\alpha_1(\cdot)$ and an
$\mathbb{H}$-predictable random measure $\alpha_2(\cdot,\cdot)$, if they exist, such that

\begin{align}\label{eq4.2}
B(t)&=\hat{B}(t) +\int_0^t\alpha_1(s)ds,\\
\tilde{N}(t,d\zeta)& = M_2(t,d\zeta) + \alpha_2(t,d\zeta) \label{eq4.2a},
\end{align}
where  $\hat{B}(t)$ is a Brownian motion with respect to $\mathbb{H}$ and $M_2(t,d\zeta)$ is an $\mathbb{H}$-martingale.
\end{problem}

 We approach this problem by studying a stochastic control problem, described as follows:

Consider a controlled jump diffusion $X(t)=X^u(t)$ of the form
\begin{align}\label{eq2.1a}
\begin{cases}
dX(t)=X(t^-)u(t)[b(t)dt+\sigma(t)dB(t)+\int_{\mathbb{R}} \gamma(t,\zeta) \tilde{N}(dt,d\zeta)]; \quad 0\leq t\leq T\\
X(0)=1,
\end{cases}
\end{align}
where $b(t), \sigma(t)$ and $\gamma(t,\zeta)$ are given $\mathbb{F}-$predictable processes.
Here $u$ is our control process, which is allowed to be $\mathbb{H}$-predictable. We say that $u$ is admissible if, in addition, the corresponding equation \eqref{eq2.1a} has a unique solution $X^{u}$ such that
\begin{equation}
X^{u}(t^-)u(t)\gamma(t,\zeta) > -1.
\end{equation}
This condition prevents $X(t)$ from jumping to a negative value.
Let $\mathcal{A}_{\mathbb{H}}$ denote this family of admissible controls $u$. The problem we study, is the following:

\begin{problem}
Find $u^{\ast}\in\mathcal{A}_{\mathbb{H}}$ such that
\begin{equation}\label{eq2.3}
\sup_{u\in\mathcal{A}_{\mathbb{H}}}\mathbb{E}[\ln(X^u(T))]=\mathbb{E}[\ln(X^{u^{\ast}}(T))]
\end{equation}
\end{problem}

\begin{remark}
A reader who is familiar with mathematical finance, will see that Problem 1.2 can be seen to represent an optimal insider portfolio problem in a financial market consisting of the following two investment possibilities:
\begin{itemize}
\item
A risk free asset with unit price $S_0(t)=1$ for all $t$\\
\item
A risky asset, with unit price $S(t)$ at time $t$ given by the SDE
\begin{align}
&dS(t)=S(t^-)[b(t)dt+\sigma(t)dB(t)+\int_{\mathbb{R}} \gamma(t,\zeta) \tilde{N}(dt,d\zeta)]; \quad 0\leq t\leq T\nonumber\\
&S(0) > 0.
\end{align}
\end{itemize}
Note, however, that we make no such interpretation here; we study Problem 1.2 purely as a stochastic control problem. Therefore questions regarding arbitrages in such a financial market are not relevant.
\end{remark}

We will solve Problem 1.2 in $3$ different ways:
\begin{itemize}
\item
By using semimartingale calculus within the filtration $\mathbb{H}$,
\item
by compensators and Hida-Malliavin calculus,
\item
by white noise calculus and the Donsker delta functional approach.
\end{itemize}

By combining the solutions by the 3 different methods we get explicit expressions for the information drift in terms of the conditional Donsker delta functional and its Hida-Malliavin derivative. See Theorem 3.1 and Theorem 4.1. \\

The results obtained in this paper are related to results presented in many earlier papers. We mention in particular \cite{I},  \cite{JY}, \cite{P}, \cite{A}, \cite{AIS}, \cite{J}, \cite{MY}, \cite{JP} and \cite{AZ}. But our method is different from the methods used in these papers, and it gives in some cases more explicit representations. Perhaps the paper which has results closest to ours, is \cite{A}. For example, Theorem 2.4 in \cite{A} corresponds to our Theorem 3.1 equation (3.15), and Theorem 2.6 in \cite{A} corresponds to our Theorem 3.1 equation (3.13).  Note, however, that our Theorem 3.1 equation (3.15) is more explicit than Theorem 2.4 in \cite{A}, because Theorem 2.4 in \cite{A} is based on the representation property (3), which is an existence result and not explicit. In Theorem 4.3 of \cite{A} a more explicit representation than in Theorem 2.4 is given, but it is based on the assumption that the difference trace exists, which need not be the case in general. Since we work in $(\mathcal{S})^*$, we do not need such an existence assumption of the corresponding Hida-Malliavin derivative.\\\\
We refer to \cite{OS} for an introduction to stochastic calculus and control of It\^ o-L\'{e}vy processes. \\

\section{The Donsker delta functional}
We first present a short introduction to the Donsker delta functional. See e.g. \cite{DO1} for more details.
\begin{definition}
Let $Y :\Omega\rightarrow\mathbb{R}$ be a random variable which also belongs to the Hida space $(\mathcal{S})^{\ast}$ of stochastic distributions. Then a continuous functional
\begin{equation}\label{donsker}
    \delta_Y(.): \mathbb{R}\rightarrow (\mathcal{S})^{\ast}
\end{equation}
is called a Donsker delta functional of $Y$ if it has the property that
\begin{equation}\label{donsker property }
    \int_{\mathbb{R}}g(y)\delta_Y(y)dy = g(Y) \quad a.s.
\end{equation}
for all (measurable) $g : \mathbb{R} \rightarrow \mathbb{R}$ such that the integral converges.
\end{definition}

The Donsker delta functional $\delta_{Y}(y)$ of a given real random variable $Y$ is related to the \emph{regular conditional distribution} of $Y$ with respect to the $\sigma$-algebra  $\mathcal{F}_t$, denoted by $Q_t(dy)=Q_t(\omega,dy)$, which is defined by the following properties (see e.g. \cite{P}):
\begin{itemize}
\item
For any Borel set $\Lambda \subseteq \mathbb{R}, Q_t(\cdot, \Lambda)$ is a version of $\mathbb{E}[\mathbf{1}_{Y \in \Lambda} | \mathcal{F}_t]$.
\item
For each fixed $\omega, Q_t(\omega,dy)$ is a probability measure on the Borel subsets of $\mathbb{R}$.
\end{itemize}

It is known that a regular conditional distribution always exists. See e. g. \cite{B}, page 79.  From the required properties of $Q_t(\omega,dy)$ it follows that
\begin{equation}
\int_{\mathbb{R}} g(y) Q_t(\omega,dy) = \mathbb{E}[g(Y) | \mathcal{F}_t];  \text{  for all bounded measurable functions  } g.
\end{equation}
If we compare this with the definition of the Donsker delta functional, we obtain the following representation of the regular conditional distribution:

\begin{proposition}
Suppose $Q_t(\omega,dy)$ is absolutely continuous with respect to Lebesgue measure on $\mathbb{R}$. Then the Donsker delta functional of $Y$,  $\delta_Y(y)$, exists in $(\mathcal{S})^*$ and we have
\begin{equation}
\frac{Q_t(\omega,dy)}{dy} = \mathbb{E}[ \delta_Y(y) | \mathcal{F}_t].
\end{equation}
\end{proposition}

The advantages with working with the Donsker delta functional, rather than the regular conditional distribution, include the following (see e.g. \cite{DO1} for details and examples):
\begin{itemize}
\item
Using white noise theory and Wick calculus one can obtain explicit formulas for the Donsker delta functional as an element of the Hida stochastic distribution space $(\mathcal{S})^*$.
\item
The Malliavin derivative has a natural extension to the \emph{Hida-Malliavin derivative} on $(\mathcal{S})^*$, and combining this extension with white noise theory and Wick calculus one can compute the Hida-Malliavin derivative of $\delta_{Y}(y)$ as an element of $(\mathcal{S})^*$.
\item
Taking condition expectation typically brings us back to $L^2(\mathbf{P})$. See Example 3.1, Example 4.1 and Example 4.2 with corollaries.
\end{itemize}
\section{The stochastic control problem}


We now turn to Problem 1.2, which we will solve by 3 different methods:

\subsection{Method 1: Enlargement of filtration.}

As pointed out in the Introduction we know by the Jacod condition  that $B(t)$ and $\tilde{N}(t,d\zeta)$ are semimartingales with respect to $\mathbb{H}$.
Then by \eqref{eq4.2} and \eqref{eq4.2a}, we can write \eqref{eq2.1a} as
\begin{align}\label{eq4.5}
\begin{cases}
dX(t)=X(t^-)u(t)[\beta(t)dt+\sigma(t)d\hat{B}(t)+\int_{\mathbb{R}}\gamma(t,\zeta)M_2(dt,d\zeta)]; \quad 0\leq t\leq T\\
X(0)=1,
\end{cases}
\end{align}
where
\begin{equation}
\beta(t)=b(t)+\sigma(t)\alpha_1(t) +\int_{\mathbb{R}}\gamma(t,\zeta)\alpha_2(t,d\zeta).
\end{equation}
This is a well-defined SDE in the semimartingale context of the $\mathbb{H}$-filtration, and we can apply classical
semimartingale calculus to solve the problem, as follows:

Let
$$\nu_{\mathbb{H}}(dt,d\zeta)$$
denote the $\mathbb{H}$-compensator of $M_2$.

Then by the It\^ o formula for semimartingales we get that the solution of the SDE \eqref{eq4.5} is
\begin{align}
X(t)&= \exp \Big( \int_0^t \{ u(s) \beta(s)-\frac{1}{2}u^2(s)\sigma^2(s) \} ds+\int_0^t\int_{\mathbb{R}}[\ln (1+u(s)\gamma(s,\zeta))-u(s)\gamma(s,\zeta)]\nu_{\mathbb{H}}(ds,d\zeta) \nonumber\\
&+\int_0^t u(s)\sigma(s) d\hat{B}(s)+\int_0^t\int_{\mathbb{R}} \ln(1+u(s)\gamma(s,\zeta)) M_2(ds,d\zeta)\Big).
\end{align}
From this we deduce that
\begin{equation}
\mathbb{E}[\ln X(T)] = \mathbb{E}[ \int_0^T \{ u(s)\beta(s)-\frac{1}{2}u^2(s)\sigma^2(s)\} ds +\int_0^T \int_{\mathbb{R}}[\ln (1+u(s)\gamma(s,\zeta))-u(s)\gamma(s,\zeta)]\nu_{\mathbb{H}}(ds,d\zeta) ].
\end{equation}

Maximizing this integrand with respect to $u(s)$ for each $s$, we get the following first order equation for the optimal $u=u^*$:
\begin{equation}\nonumber
\{ \beta(s)-u^*(s)\sigma^2(s)\}ds+ \int_{\mathbb{R}}[\frac{\gamma(s,\zeta)}{1+u^*(s)\gamma(s,\zeta)}-\gamma(s,\zeta)]\nu_{\mathbb{H}}(ds,d\zeta)=0,
\end{equation}
i.e.
\begin{equation}\label{eq4.9}
\{ b(s)-u^*(s)\sigma^2(s)+\sigma(s)\alpha_1(s)\}ds- \int_{\mathbb{R}}\frac{u^*(s)\gamma^2(s,\zeta)}{1+u^*(s)\gamma(s,\zeta)}\nu_{\mathbb{H}}(ds,d\zeta)=-\int_{\mathbb{R}}\gamma(s,\zeta)\alpha_2(ds,d\zeta),
\end{equation}
provided that
\begin{equation}\label{eq4.9a}
1+u^*(s)\gamma(s,\zeta) \neq 0.
\end{equation}

\subsection{Method 2: Using Hida-Malliavin calculus}
 Problem 2.1 is also studied in \cite{DMOP}. Using Hida-Malliavin calculus, combined with uniqueness of semimartingale decompositions, it follows by Theorem 16 in \cite{DMOP} that
if an optimal portfolio $u^*$ exists, then it satisfies the equation
\begin{equation}\label{eq4.10}
\{ b(s)-u^*(s)\sigma^2(s)+\sigma(s)\alpha_1(s)\}ds -\int_{\mathbb{R}}\frac{u^*(s)\gamma^2(s,\zeta)}{1+u^*(s)\gamma(s,\zeta)}\nu(d\zeta) ds=\int_{\mathbb{R}}\frac{\gamma(s,\zeta)}{1+u^*(s)\gamma(s,\zeta)} (\nu(d\zeta)ds-\nu_{\mathbb{H}}(ds,d\zeta)),
\end{equation}
provided that
\begin{equation}\label{eq4.11}
1+u^*(s)\gamma(s,\zeta) \neq 0.
\end{equation}

\subsection{Method 3: Using white noise theory and the Donsker delta functional.}

By Theorem 6.8 in \cite{DO1} it follows that the optimal control $u^{\ast}$ is a solution of:
 \begin{align}\label{eq4.13}
 b(s)-u^*(s)\sigma^2(s)+\sigma(s)\Phi(s)- \int_{\mathbb{R}} \frac{u^{\ast}(s) \gamma^2(s,\zeta)}{1+u^{\ast}(s)\gamma(s,\zeta)} \nu(d\zeta)
=- \int_{\mathbb{R}}\frac{\gamma(s,\zeta)}{1+u^{\ast}(s)\gamma(s,\zeta)}\Psi(s,\zeta)\nu(d\zeta),
\end{align}
where
\begin{align}\label{eq4.14}
\Phi(s):=\frac{\mathbb{E}[D_{s}\delta_Y(y)|\mathcal{F}_s]_{y=Y}}{\mathbb{E}[\delta_Y(y)|\mathcal{F}_s]_{y=Y}}\quad \text{ and }\quad
\Psi(s,\zeta):= \frac{\mathbb{E}[D_{s,\zeta}\delta_Y(y)|\mathcal{F}_s]_{y=Y}}{\mathbb{E}[\delta_Y(y)|\mathcal{F}_s]_{y=Y}}.
\end{align}
\vskip 0.5cm

Combining the results from the 3 methods above we get the following theorem. Note that equation \eqref{eq4.15}  may be regarded as a Donsker delta analogue of Theorem 2.6 in \cite{A}, and equation \eqref{eq4.18a} a Donsker delta analogue of Theorem 2.4 in \cite{A}. However, as explained in the Introduction, our expressions are more explicit than the analogue results in \cite{A}.

\begin{theorem}[$\mathbb{H}$ compensator and $\mathbb{H}$-semimartingale decomposition]
In the \\
$\mathbb{H}-$semimartingale decomposition \eqref{eq4.2},\eqref{eq4.2a} the following hold:\\
\begin{itemize}
\item
(i) The $\mathbb{H}$-compensator of $\tilde{N}$ is given by
\begin{equation}\label{eq4.15}
\nu_{\mathbb{H}}(ds,d\zeta)=(1+\Psi(s,\zeta))\nu(d\zeta)ds.
\end{equation}
\item
(ii) The process $\alpha_1(t)$ in the $\mathbb{H}$-semimartingale decomposition \eqref{eq4.2} of $B$ is
\begin{equation}\label{eq3.14a}
\alpha_1(t)=\Phi(t)
\end{equation}
\item
(iii) The process $\alpha_2(t,\zeta)$ in the $\mathbb{H}$-semimartingale decomposition \eqref{eq4.2a} of $\tilde{N}$ is
\begin{equation}\label{eq4.18a}
\alpha_2(t,d\zeta)=\Psi(t,\zeta)\nu(d\zeta).
\end{equation}
\end{itemize}
\end{theorem}

\dproof\\
\emph{(i) and (ii)}: Subtracting \eqref{eq4.10} from \eqref{eq4.13} we get:
\begin{equation}\label{eq4.15a}
\sigma(s)(\Phi(s)-\alpha_1(s))ds= \int_{\mathbb{R}}\frac{\gamma(s,\zeta)}{1+u^{\ast}(s)\gamma(s,\zeta)}\nu_{\mathbb{H}}(ds,d\zeta)-\int_{\mathbb{R}}\frac{\gamma(s,\zeta)}{1+u^{\ast}(s)\gamma(s,\zeta)}(1+\Psi(s,\zeta))\nu(d\zeta)ds.
\end{equation}
Since this holds for all $\gamma(s,\zeta)$, we conclude that \eqref{eq4.15} holds. Moreover, we can conclude that
\eqref{eq3.14a} holds also.\\

\noindent
\emph{(iii)}: Substituting \eqref{eq4.15} into \eqref{eq4.9} we get
\begin{align}\label{eq4.16}
&b(s)-u^*(s)\sigma^2(s)+\sigma(s)\alpha_1(s)-\int_{\mathbb{R}}\frac{u^*(s)\gamma^2(s,\zeta)}{1+u^*(s)\gamma(s,\zeta)}\nu(d\zeta)\nonumber\\
&-\int_{\mathbb{R}}\frac{u^*(s)\gamma^2(s,\zeta)}{1+u^*(s)\gamma(s,\zeta)}\Psi(s,\zeta)\nu(d\zeta)=-\int_{\mathbb{R}}\gamma(s,\zeta)\alpha_2(s,d\zeta).
\end{align}
Substituting  \eqref{eq4.13} into \eqref{eq4.16} we get
\begin{align}
&\int_{\mathbb{R}}\frac{\gamma(s,\zeta)}{1+u^*(s)\gamma(s,\zeta)}\Psi(s,\zeta)\nu(d\zeta)\nonumber\\
&+\int_{\mathbb{R}}\frac{u^*(s)\gamma^2(s,\zeta)}{1+u^*(s)\gamma(s,\zeta)}\Psi(s,\zeta)\nu(d\zeta)=\int_{\mathbb{R}}\gamma(s,\zeta)\alpha_2(s,d\zeta)
\end{align}
i.e.
\begin{align}
\int_{\mathbb{R}}\gamma(s,\zeta)\Psi(s,\zeta)\nu(d\zeta)=\int_{\mathbb{R}}\gamma(s,\zeta)\alpha_2(s,d\zeta).
\end{align}

Since this holds for all $\gamma(t,\zeta)$ and $\alpha_2$ does not depend on $\gamma$, we deduce that \eqref{eq4.18a} holds.
\fproof

\example{[Semimartingale decomposition with respect to $\mathbb{H}=\mathbb{H}_Y$, when $Y=Y(T_0)$ and $Y$ is an It\^ o-L\' evy process]}\\

We refer to \cite{DO1} and the references therein for more details in this example. \\
Consider the case when $Y$ is a first order chaos random variable of the form

\begin{equation}\label{eq4.22}
    Y = Y (T_0); \text{ where } Y (t) =\int_0^t\sigma(s)dB(s)+\int_0^t\int_{\mathbb{R}}\theta(s,\zeta)\tilde{N}(ds,d\zeta), \mbox{ for } t\in [0,T_0]
\end{equation}
for some deterministic functions $\sigma\neq0, \theta$ satisfying
\begin{equation}\label{}
    \int_0^{T_0} \{ \sigma^2(t)+\int_{\mathbb{R}}\theta^2(t,\zeta)\nu(d\zeta)\} dt<\infty \text{ a.s. }
\end{equation}

In this case the Donsker delta functional of $Y$ exists in $(\mathcal{S})^{\ast}$ and is given
by

\begin{eqnarray}\label{eq2.7}
   \delta_Y(y)&=&\frac{1}{2\pi}\int_{\mathbb{R}}\exp^{\diamond}\big[ \int_0^{T_0}\int_{\mathbb{R}}(e^{ix\theta(s,\zeta)}-1)\tilde{N}(ds,d\zeta)+ \int_0^{T_0}ix\sigma(s)dB(s)  \nonumber\\
   &+&  \int_0^{T_0}\{\int_{\mathbb{R}}(e^{ix\theta(s,\zeta)}-1-ix\theta(s,\zeta))\nu(d\zeta)-\frac{1}{2}x^2\sigma^2(s)\}ds-ixy\big]dx.
\end{eqnarray}

From this one can deduce that
\begin{align}
\mathbb{E}[\delta_Y(y)|\mathcal{F}_t] &=\frac{1}{2\pi}\int_{\mathbb{R}}\exp\big[\int_0^t\int_{\mathbb{R}}ix\theta(s,\zeta)\tilde{N}(ds,d\zeta) +\int_0^t ix\sigma(s)dB(s)\nonumber\\
   &+\int_t^{T_0}\int_{\mathbb{R}}(e^{ix\theta(s,\zeta)}-1-ix\theta(s,\zeta))\nu(d\zeta)ds-\int_t^{T_0}\frac{1}{2}x^2\sigma^2(s)ds-ixy\big]dx.
   \end{align}
and
\begin{align}
&\mathbb{E}[D_{t,\zeta}\delta_Y(y)|\mathcal{F}_t]\nonumber\\
 &=\frac{1}{2\pi}\int_{\mathbb{R}}\exp\big[\int_0^t\int_{\mathbb{R}}ix\theta(s,\zeta)\tilde{N}(ds,d\zeta) +\int_0^t ix\sigma(s)dB(s)\nonumber\\
   &+\int_t^{T_0}\int_{\mathbb{R}}(e^{ix\theta(s,\zeta)}-1-ix\theta(s,\zeta))\nu(d\zeta)ds-\int_t^{T_0}\frac{1}{2}x^2\sigma^2(s)ds-ixy\big]\nonumber\\
   &\times(e^{ix\theta(t,z)}-1)dx.
  \end{align}

and
\begin{align}
&\mathbb{E}[D_{t}\delta_Y(y)|\mathcal{F}_t]\nonumber\\
&= \frac{1}{2\pi}\int_{\mathbb{R}}\exp\big[\int_0^t\int_{\mathbb{R}}ix\theta(s,\zeta)\tilde{N}(ds,d\zeta) +\int_0^t ix\sigma(s)dB(s)\nonumber\\
   &+\int_t^{T_0}\int_{\mathbb{R}}(e^{ix\theta(s,\zeta)}-1-ix\theta(s,\zeta))\nu(d\zeta)ds-\int_t^{T_0}\frac{1}{2}x^2\sigma^2(s)ds-ixy\big]ix\sigma(t)dx.
\end{align}

  Therefore we get the following result, which may be viewed as an explicit Donsker delta functional version of Theorem 3.5 in \cite{AZ}:

  \begin{theorem}
  If $Y=Y(T_0)$ as in \eqref{eq4.22}, then the process $\alpha_2$ in the $\mathbb{H}$-semimartingale decomposition \eqref{eq4.2a} of $\tilde{N}$ is
  \begin{align}\label{eq4.26}
  \alpha_2(t,d\zeta)=\Psi(t,\zeta) \nu(d\zeta):= \frac{\mathbb{E}[D_{s,\zeta}\delta_Y(y)|\mathcal{F}_s]_{y=Y}}{\mathbb{E}[\delta_Y(y)|\mathcal{F}_s]_{y=Y}}\nu(d\zeta)=\frac{\int_{\mathbb{R}} F(t,x,Y) (e^{ix\theta(t,\zeta)}-1)dx}{\int_{\mathbb{R}} F(t,x,Y) dx} \nu(d\zeta)
\end{align}
where
\begin{align}
F(t,x,y)&=\exp\big[\int_0^t\int_{\mathbb{R}}ix\theta(s,\zeta)\tilde{N}(ds,d\zeta) +\int_0^t ix\sigma(s)dB(s)\nonumber\\
&+\int_t^{T_0}\int_{\mathbb{R}}(e^{ix\theta(s,\zeta)}-1-ix\theta(s,\zeta))\nu(d\zeta)ds-\int_t^{T_0}\frac{1}{2}x^2\sigma^2(s)ds-ixy\big].
\end{align}
   \end{theorem}

\section{ Some special cases}
\subsection{The $\mathbb{H}$-semimartingale decomposition of Brownian motion}

From Theorem 3.1 we have
\begin{equation}\label{eq2.15}
\alpha_1(t)=\frac{\mathbb{E}[D_t\delta_Y(y)|\mathcal{F}_t]_{y=Y}}{\mathbb{E}[\delta_Y(y)|\mathcal{F}_t]_{y=Y}}.
\end{equation}
The right hand side of (\ref{eq2.15}) can in often be computed explicitly  using known expressions of the Donsker delta functional and its Hida-Malliavin derivative given in \cite{DO1}. Here is an example, which is well-known, but added as an illustration of our method:
\begin{example}
Suppose
\begin{equation}\label{eq5.47}
    Y = Y (T_0); \text{ where } Y (t) =\int_0^t\beta(s)dB(s), \mbox{ for } t\in [0,T_0]
\end{equation}
for some deterministic function $\beta\in \mathbf{L}^2[0,T_0]$ with
\begin{equation}\label{}
    \|\beta\|^2_{[t,T]} :=\int_t^T\beta(s)^2ds>0 \mbox{ for all } t\in[0,T].
\end{equation}
In this case we know that the Donsker delta functional is given
by
\begin{equation}\label{}
    \delta_{Y}(y)=(2\pi v)^{-\frac{1}{2}}\exp^{\diamond}[-\frac{(Y-y)^{\diamond2}}{2v}]
\end{equation}
where we have put $v :=\|\beta\|^2_{[0,T_0]}$. See e.g. \cite{AaOU}, Proposition $3.2$.
Using the Wick rule when taking conditional expectation, using the martingale property of
the process $Y (t)$ and applying Lemma $3.7$ in \cite{AaOU} we get
\begin{eqnarray}\label{eq5.50}
   \mathbb{E}[\delta_Y(y)|\mathcal{F}_t]&=&(2\pi v)^{-\frac{1}{2}}\exp^{\diamond}[-\mathbb{E}[\frac{(Y(T_0)-y)^{\diamond 2}}{2v}|\mathcal{F}_t]] \nonumber \\
   &=& (2\pi \|\beta\|^2_{[0,T_0]})^{-\frac{1}{2}}\exp^{\diamond}[- \frac{(Y(t)-y)^{\diamond 2}}{2\|\beta\|^2_{[0,T_0]}}]\nonumber \\
   &=& (2\pi \|\beta\|^2_{[t,T_0]})^{-\frac{1}{2}} \exp[- \frac{(Y(t)-y)^2}{2\|\beta\|^2_{[t,T_0]}}].
\end{eqnarray}
Similarly, by the Wick chain rule and Lemma $3.8$ in \cite{AaOU} we get, for $t \in [0,T],$
\begin{eqnarray}\label{eq5.51}
  \mathbb{E}[D_t\delta_Y(y)|\mathcal{F}_t] &=&-\mathbb{E}[(2\pi v)^{-\frac{1}{2}}\exp^{\diamond}[- \frac{(Y(T_0)-y)^{\diamond 2}}{2v}]\diamond\frac{Y(T_0)-y}{v}\beta(t)|\mathcal{F}_t] \nonumber\\
   &=&-(2\pi v)^{-\frac{1}{2}} \exp^{\diamond}[- \frac{(Y(t)-y)^{\diamond 2}}{2v}]\diamond\frac{Y(t)-y}{v}\beta(t)\nonumber \\
   &=& -(2\pi \|\beta\|^2_{[t,T_0]})^{-\frac{1}{2}}\exp[- \frac{(Y(t)-y)^2}{2\|\beta\|^2_{[t,T_0]}}]\frac{Y(t)-y}{\|\beta\|^2_{[t,T_0]}}\beta(t).
\end{eqnarray}
We conclude that in this case we have, by Theorem 3.1,
\begin{equation}
\alpha(s)=\Phi(s) = \frac{Y(T_0)-Y(s)}{\|\beta\|^2_{[s,T_0]}}\beta(s),\quad s \in [0,T].
\end{equation}

\end{example}

\subsection{The $\mathbb{H}$-semimartingale decomposition of a Poisson process}
Next, consider the special case when $N$ is a Poisson process with constant intensity $\lambda>0$, so that  $\tilde{N}(t) := N(t)-\lambda t$.
By Theorem 3.1 we then get
 \begin{theorem}{[The $\mathbb{H}$-enlargement of filtration formula for a Poisson process]}
 \begin{itemize}
 \item
 The unique $\mathbb{H}$-predictable compensator $\rho$ of the compensated Poisson process $\tilde{N}(t)=N(t)-\lambda t$ is given by
 \begin{equation}\label{eq3.20}
 \rho(t)=\lambda(1+\Psi(t,1)).
 \end{equation}
\item
The $\mathbb{H}$-semimartingale decomposition of $\tilde{N}$ is
\begin{equation}\label{eq3.21}
\tilde{N}(t)=M(t)+\int_0^t\alpha(s) ds,
\end{equation}
where $M(t)$ is an $\mathbb{H}$-martingale, and
\begin{equation}\label{eq3.22}
\alpha(s)=\lambda\Psi(s,1), \quad 0 \leq s \leq T.
 \end{equation}
 \end{itemize}
\end{theorem}

\begin{example}
We refer to \cite{DO1} and the references therein for more details in the following.\\
Assume that the inside information from time $t=0$ is the value of $Y=Y(T_0)$, with
\begin{equation}\label{eq3.24}
 Y(t)=\theta B(t)+\tilde{N}(t); \quad 0 \leq t \leq T_0
 \end{equation}
where $\theta \neq 0$ is a constant.
In this case the
L\'{e}vy measure is $\nu(d\zeta) = \lambda\delta_1(d\zeta)$ since the jumps are of size $1$.
Thus we have, with $i=\sqrt{-1}$,
\begin{equation}\label{eq2.20a}
   \delta_Y(y)=\frac{1}{2\pi}\int_{\mathbb{R}}\exp^{\diamond}\big[(e^{ix}-1)\tilde{N}(T_0)+i x \theta B(T_0)+\lambda T_0 (e^{ix}-1-ix)-\frac{1}{2}x^2 \theta^2 T_0-ixy\big]dx
\end{equation}

From this we get:
\begin{equation}\label{eq2.21a}
    \mathbb{E}[\delta_Y(y)|\mathcal{F}_t]=\frac{1}{2\pi}\int_{\mathbb{R}}\exp\big[ ix\tilde{N}(t)+i x \theta B(t)+\lambda(T_0-t)(e^{ix}-1-ix)-\frac{1}{2}x^2 \theta^2 (T_0-t)-ixy\big]dx
\end{equation}

\begin{equation}\label{eq2.22a}
    \mathbb{E}[D_{t,1}\delta_Y(y)|\mathcal{F}_t]=\frac{1}{2\pi}\int_{\mathbb{R}}\exp\big[ ix\tilde{N}(t)+i x \theta B(t)+\lambda(T_0-t)(e^{ix}-1-ix)-\frac{1}{2}x^2 \theta^2 (T_0-t)-ixy\big](e^{ix}-1)dx.
\end{equation}
and
\begin{equation}\label{eq2.22a}
\mathbb{E}[D_{t}\delta_Y(y)|\mathcal{F}_t]=i\frac{1}{2\pi}\theta\int_{\mathbb{R}}\exp\big[ ix\tilde{N}(t)+i x \theta B(t)+\lambda(T_0-t)(e^{ix}-1-ix)-\frac{1}{2}x^2 \theta^2 (T_0-t)-ixy\big]xdx.
\end{equation}

By \eqref{eq4.14} we conclude that
\begin{equation}
\Psi(t,1)= \frac{\int_{\mathbb{R}}\exp\big[ ix\tilde{N}(t)+i x \theta B(t)+\lambda(T_0-t)(e^{ix}-1-ix)-\frac{1}{2}x^2 \theta^2 (T_0-t)-ixY\big](e^{ix}-1)dx}{\int_{\mathbb{R}}\exp\big[ ix\tilde{N}(t)+i x \theta B(t)+\lambda(T_0-t)(e^{ix}-1-ix)-\frac{1}{2}x^2 \theta^2 (T_0-t)-ixY\big]dx}
\end{equation}
Note that if we put
\begin{equation}
Z(t,x) :=\exp\big[ ix\tilde{N}(t)+i x \theta B(t)+\lambda(T_0-t)(e^{ix}-1-ix)-\frac{1}{2}x^2 \theta^2 (T_0-t)-ixY\big],
\end{equation}
then
\begin{equation}
\frac{\partial Z}{\partial x}(t,x)= i Z(t,x) [\tilde{N}(t)+\theta B(t)+ix \theta^2 (T_0-t)+\lambda (T_0-t) (e^{ix}-1) - Y],
\end{equation}
or
\begin{equation}
Z(t,x)(e^{ix}-1)= \frac{-i\frac{\partial Z}{\partial x} - Z(t,x) [\tilde{N}(t)+\theta B(t)+ix \theta^2 (T_0-t) -Y]}{\lambda(T_0-t)}
\end{equation}\\

Therefore, since $Z(t,x) = 0$ at $x=\infty$ and at $x= -\infty$, we get by (3.27) the following result:

\begin{proposition}
Suppose $Y =Y(T_0)$, with $Y(t)$ as in \eqref{eq3.24}. Then we have

\begin{align}\label{eq3.33}
\Psi(t,1)&:=\frac {\mathbb{E}[D_{t,1}\delta_Y(y)|\mathcal{F}_t]_{y=Y}}{\mathbb{E}[\delta_Y(y)|\mathcal{F}_t]_{y=Y}}=\frac{Y-\theta B(t)-\tilde{N}(t)}{\lambda(T_0-t)}-i\theta^2\frac{\int_{\mathbb{R}} xZ(t,x)dx}{\lambda\int_{\mathbb{R}} Z(t,x)dx}\nonumber\\
&=\frac{Y-\theta B(t)-\tilde{N}(t)}{\lambda(T_0-t)}-\frac{\theta^2}{\lambda}\Phi(t),
\end{align}
where
\begin{align}
\Phi (t)&:= \frac{\mathbb{E}[D_t \delta_Y(y)|\mathcal{F}_t]_{y=Y}}{\mathbb{E}[\delta_Y(y)|\mathcal{F}_t]_{y=Y}} \nonumber\\
&= \frac{\int_{\mathbb{R}}\exp\big[ ix\tilde{N}(t)+i x \theta B(t)+\lambda(T_0-t)(e^{ix}-1-ix)-\frac{1}{2}x^2 \theta^2 (T_0-t)-ixY\big]xdx}{\int_{\mathbb{R}}\exp\big[ ix\tilde{N}(t)+i x \theta B(t)+\lambda(T_0-t)(e^{ix}-1-ix)-\frac{1}{2}x^2 \theta^2 (T_0-t)-ixY\big]dx}.
\end{align}
Hence, by \eqref{eq3.22} we then get
\begin{equation}
\alpha(t)=\frac{Y-\theta B(t)-\tilde{N}(t)}{T_0-t}-\theta^2\Phi(t).
\end{equation}
\end{proposition}

In particular, by letting $\theta \rightarrow 0$ in \eqref{eq3.33} we get as a special case the result from Example 3.1:

\begin{coro}
Suppose $Y=\tilde{N}(T_0)$. Then
\begin{equation}
\Psi(t,1)=\frac{Y-\tilde{N}(t)}{\lambda(T_0-t)},
\end{equation}
which by \eqref{eq3.22} gives
\begin{equation}
\alpha(t)=\frac{\tilde{N}(T_0) - \tilde{N}(t)}{T_0-t}.
\end{equation}
\end{coro}

\end{example}

\vskip 0.3cm
\noindent \textbf{Acknowledgment}\\
We are grateful to Monique Jeanblanc for helpful comments.

\end{document}